
\input amssym.def
\input amssym
\input psfig
\magnification=1100
\baselineskip = 0.25truein
\lineskiplimit = 0.01truein
\lineskip = 0.01truein
\vsize = 8.5truein
\voffset = 0.2truein
\parskip = 0.10truein
\parindent = 0.3truein
\settabs 12 \columns
\hsize = 5.4truein
\hoffset = 0.4truein

\setbox\strutbox=\hbox{%
\vrule height .708\baselineskip
depth .292\baselineskip
width 0pt}
\font\caps=cmcsc10

\font\bigtenrm=cmr10 at 14pt

\def\sqr#1#2{{\vcenter{\vbox{\hrule height.#2pt
\hbox{\vrule width.#2pt height#1pt \kern#1pt
\vrule width.#2pt}
\hrule height.#2pt}}}}
\def\square{\mathchoice\sqr46\sqr46\sqr{3.1}6\sqr{2.3}4}

\centerline{\bigtenrm THE HEEGAARD GENUS}
\centerline{\bigtenrm OF AMALGAMATED 3-MANIFOLDS}
\tenrm
\vskip 14pt
\centerline{MARC LACKENBY}
\vskip 18pt

\tenrm
\centerline{\caps 1. Introduction}
\vskip 6pt

When studying Haken 3-manifolds, one is led naturally
to the following construction: the amalgamation of
two 3-manifolds $M$ and $M'$ via a homeomorphism
between their boundaries. In this paper, we study
the behaviour of Heegaard genus under this operation.
We show that, provided the gluing homeomorphism
is `sufficiently complicated' and $M$ and $M'$
satisfy some standard conditions, then the Heegaard
genus of the amalgamated manifold is completely
determined by the Heegaard genus of $M$ and $M'$
and the genus of their common boundary.
Recall that a 3-manifold is {\sl simple} if it is
compact, orientable, irreducible, atoroidal, acylindrical
and has incompressible boundary. We denote the Heegaard
genus of a 3-manifold $M$ by $g(M)$.

\noindent {\bf Main Theorem.} {\sl Let $M$ and $M'$
be simple 3-manifolds, and let $h \colon
\partial M \rightarrow S$ and $h' \colon S \rightarrow
\partial M'$ be homeomorphisms with some connected surface
$S$ of genus at least two. Let $\psi \colon S \rightarrow S$ be a pseudo-Anosov
homeomorphism. Then, provided $|n|$ is sufficiently large,
$$g(M \cup_{h' \psi^n h} M') = g(M) + g(M') - g(S).$$
Furthermore, any minimal genus Heegaard splitting
for $M \cup_{h' \psi^n h} M'$ is obtained from
splittings of $M$ and $M'$ by amalgamation, and hence
is weakly reducible.}

The amalgamation of two Heegaard splittings, referred to
in the above theorem, was
defined by Schultens [10]. We recall it here. 
Since $\partial M$ and $\partial M'$ are
assumed to be connected,
the Heegaard splittings of $M$ and $M'$ divide 
each manifold into a compression body
and a handlebody. Each compression body is a copy 
of $S \times I$ with 1-handles
attached. Extend these 1-handles vertically through
$S \times I$ so that they are attached to $\partial M$
and $\partial M'$ respectively. We may assume that their
attaching discs are disjoint when the manifolds are glued.
Attach the boundaries of these 1-handles to the copy
of $S$ in $M \cup M'$, and remove the interiors of
the attaching discs. The resulting surface is a Heegaard
surface for $M \cup M'$, which is
said to be obtained from the splittings of
$M$ and $M'$ by {\sl amalgamation}. (See Figure 1.) By calculating
the genus of this surface, we obtain the
inequality
$$g(M \cup M') \leq g(M) + g(M') - g(S).$$

\vskip 12pt
\centerline{\psfig{figure=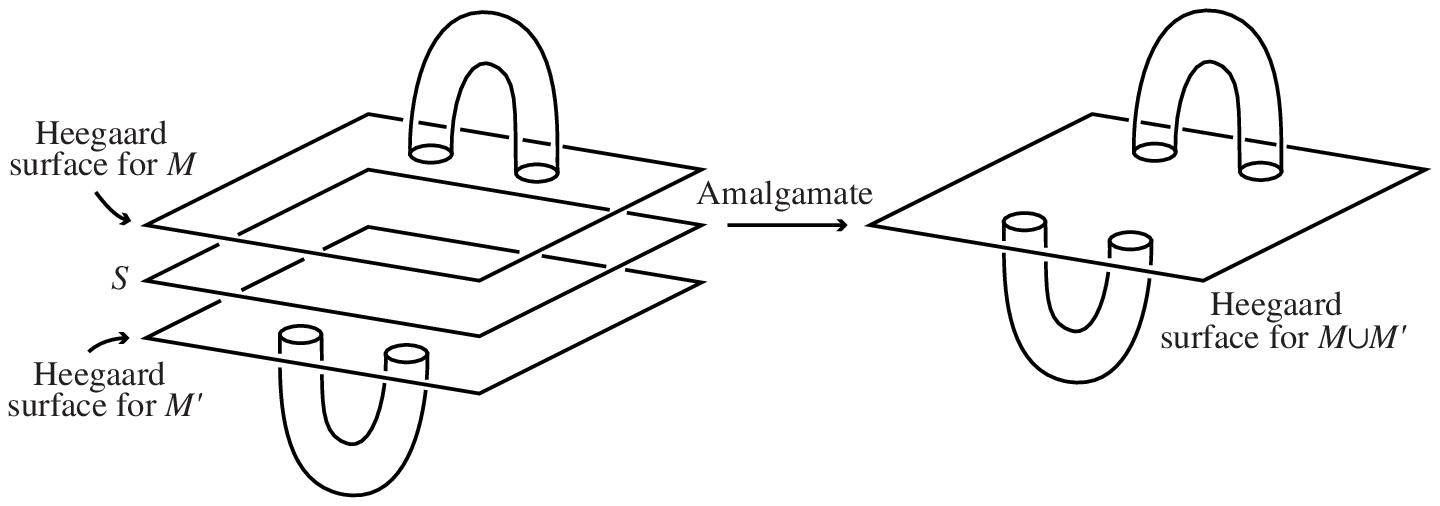,width=4.3in}}
\vskip 12pt
\centerline{Figure 1.}

Inequalities going in the other direction have been
discovered by Johannson [3] and Schultens [11]
who proved, respectively, that
$$\eqalign{
g(M \cup M') &\geq {1 \over 5} g(M) + {1 \over 5} g(M') -
{2 \over 5} g(S), \cr
g(M \cup M') &\geq {1 \over 3} g(M) + {1 \over 3} g(M') -
{4 \over 3} g(S) + {5 \over 3}. \cr}$$
Clearly, in general, one will not be able
to determine $g(M \cup M')$ precisely in terms
of $g(M)$, $g(M')$ and $g(S)$. Inequalities
such as the above will have to suffice. It is
therefore slightly surprising that, for complicated
gluings, an exact formula as in the main theorem should
hold.

\vskip 18pt
\centerline {\caps 2. Proof of the main theorem}
\vskip 6pt

The amalgamated manifold $M \cup M'$ is Haken,
atoroidal and not Seifert fibred. So, by Thurston's geometrisation theorem, 
it admits a hyperbolic structure. In [12], Soma gave a careful
analysis of its geometry. Soma proved that one can
find a point $x_n$ in each
$M \cup_{h' \psi^n h} M'$ such that the based manifolds
$(M \cup_{h' \psi^n h} M', x_n)$ converge in the
Gromov-Hausdorff topology to the infinite cyclic
cover of the hyperbolic fibred 3-manifold with monodromy $\psi$.
Furthermore, any fibre in the limit space pulls back
to a surface isotopic to the copy of $S$ in $M \cup_{h' \psi^n h} M'$,
provided $|n|$ is sufficiently large.
Hence, we deduce that, provided $|n|$ is sufficiently
large, one may find in $M \cup_{h' \psi^n h} M'$
an arbitrarily large number of parallel copies of $S$, 
such that any two adjacent copies have distance at
least one from each other. We denote the
product region in $M \cup_{h' \psi^n h} M'$ between the extreme copies
of $S$ by $S \times I$. Then we may also ensure that
there is an $\epsilon > 0$, independent of $n$,
such that $S \times I$ lies in the $\epsilon$-thick
part of $M \cup_{h' \psi^n h} M'$.

Now consider a minimal genus Heegaard surface $F$
for $M \cup M'$. Note that $g(F) \leq g(M) + g(M') - g(S)$.
From $F$, we construct
(as in [7], [8] or [4]) a generalised Heegaard splitting
$\{ F_1, \dots, F_m \}$ with the following properties:
\item{$\bullet$} $F_j$ is strongly irreducible, for each odd $j$;
\item{$\bullet$} $F_j$ is incompressible and has no
2-sphere components, for each even $j$;
\item{$\bullet$} $F_j$ and $F_{j+1}$ are not parallel
for any $j$;
\item{$\bullet$} $\sum_j (-1)^j \chi(F_j) = -\chi(F)$;
\item{$\bullet$} $|\chi(F_j)| \leq |\chi(F)|$ for each $j$.

\noindent Let $F_+ = F_1 \cup \dots \cup F_m$.
The third and fourth conditions imply that
$m \leq |\chi(F)|$, and hence the fifth gives that
$|\chi(F_+)| \leq |\chi(F)|^2$,
a bound which is independent of $n$.
One can obtain $F$ back from $F_+$
by amalgamating $F_1$ and $F_3$, then
amalgamating this with $F_5$, and so on.

By theorems of Schoen and Yau [9], Freedman, Hass and Scott [1]
and Pitts and Rubinstein [5],
each component of $F_+$ may be isotoped to a minimal
surface or to the double cover of a minimal
non-orientable surface (possibly with a small
tube attached in the case of an odd surface). 
Furthermore, after these isotopies,
any two components are either equal or disjoint.
Each complementary region of $F_+$ after the
isotopies corresponds to one before, but some
product complementary regions may have been
collapsed. In particular, each complementary
region afterwards is a compression body.

We would like to apply Proposition 6.1 of [4], which 
gives a constant $k$, such that 
each component $F'$ of $F_+$ has diameter
at most $k |\chi(F')|$. (Here, we are using
the path metric on $F_+$ arising from its induced
Riemannian metric.) However, $k$ depends on
a positive lower bound for the injectivity radius
of the ambient manifold. It is not immediately
clear from Soma's paper whether there is such a
bound that is independent of $n$. We will
therefore present a variant of Proposition 6.1 of [4].
Let $\delta$ be $2\epsilon + 1+ \epsilon/\pi$.
We claim that we can cover $F_+ \cap (S \times I)$
with regions, each of which has diameter
at most $\delta$ in $F_+$, and so that
the total number of regions is at most
$k |\chi(F_+)|$, for some constant $k$
independent of $n$. These regions will be of
two types. 

Let $(F_+)_{[\epsilon, \infty)}$
and $(F_+)_{(0, \epsilon]}$ be the
$\epsilon$-thick and $\epsilon$-thin parts
of $F_+$. Let $\Gamma$ be a maximal
collection of disjoint (not necessarily simple)
closed geodesics in $F_+$, each with length
less than $\epsilon$. The first type of
region will consist of those points within
$\epsilon/2 + \epsilon/(2 \pi) + 1/2$ of
some component of $\Gamma$. Clearly,
each such region has diameter at most
$\delta$. Claim 3 in the proof of
Proposition 6.1 of [4] gives that there
are at most $4 |\chi(F_+)|$ geodesics in
$\Gamma$, and hence at most $4 |\chi(F_+)|$ 
such regions. The argument of Claims 2 and 1 there
also gives that these regions cover
$(F_+)_{(0, \epsilon]} \cap (S \times I)$.
This uses the assumption that
$S \times I$ lies in the $\epsilon$-thick
part of $M \cup M'$.

Now pick a maximal collection of
points in $(F_+)_{[\epsilon, \infty)} \cap (S \times I)$,
no two of which are less than $\epsilon$
apart in $F_+$. Then the $\epsilon$-balls
around these points cover $(F_+)_{[\epsilon, \infty)}
\cap (S \times I)$, and there are at
most $(\cosh (\epsilon/2) - 1)^{-1} |\chi(F_+)|$ 
such balls. Letting these be the other
type of region, we have established
the claim.

We claim that one of the parallel copies of $S$
is disjoint from $F_+$, when $|n|$ is
sufficiently large. Since each of the
regions into which we have divided $F_+
\cap (S \times I)$ has uniformly bounded diameter,
there is a uniform upper bound on the
number of copies of $S$ it can intersect.
There is also a uniform upper bound on
the number of such regions. Hence, there
is a uniform upper bound on the number
of copies of $S$ that $F_+$ can intersect.
When $|n|$ is sufficiently large, there
are more copies of $S$ than this bound.
This proves the claim. (See Figure 2.)

So, some copy of $S$ lies in
the complement of $F_+$, which is
a collection of compression bodies. Since
$S$ is incompressible, it must be parallel
to a component of $F_j$ for some even $j$.
Thus, if we were to cut $M \cup M'$ along
$S$, we would obtain
generalised Heegaard splittings for $M$
and $M'$. Amalgamate each of these,
to form Heegaard surfaces $\tilde F$ and
$\tilde F'$ for $M$ and $M'$. Then,
$F$ is obtained by amalgamating
$\tilde F$ and $\tilde F'$ along $S$.

This implies that $g(F) = g(\tilde F) + g(\tilde F') - 
g(S) \geq g(M) + g(M') - g(S)$. Since we
already have the opposite inequality, the
theorem is proved. $\square$

\vskip 12pt
\centerline{\psfig{figure=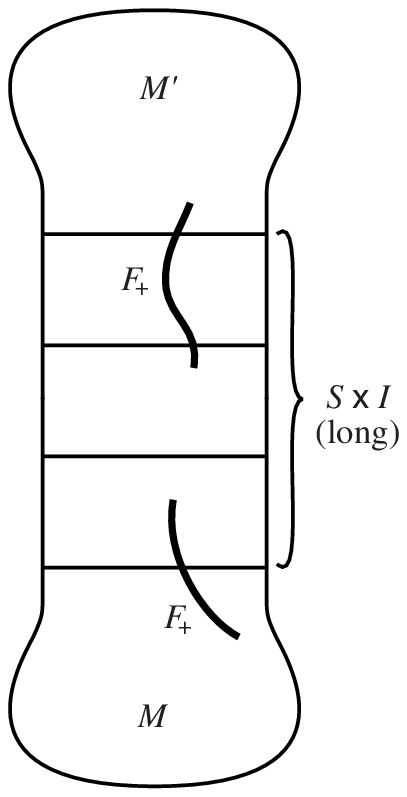,width=1.5in}}
\vskip 18pt
\centerline{Figure 2.}

\vskip 12pt
\centerline{\caps 3. Generalisations}
\vskip 6pt

The main theorem is not the most general possible
statement one can make. In fact, the proof gives the
following stronger result.

\noindent {\bf Theorem.} {\sl Let $M$, $M'$, $S$,
$h$, $h'$ and $\psi$ be as in the main theorem. Then for each
$g >0$, there is an $N > 0$ with the following
property: if $|n| \geq N$, then any genus
$g$ splitting for $M \cup_{h' \psi^n h} M'$
is obtained from splittings of $M$ and $M'$
by amalgamation. In particular, it is weakly
reducible.}

There is a related way of building a Haken
3-manifold via gluing: one can start with a single
simple 3-manifold $M$, and glue two of its boundary components
via an orientation-reversing homeomorphism. 
In this case, we obtain a similar result
to the main theorem, but do not obtain
a precise equality.

\noindent {\bf Theorem.} {\sl Let $M$
be a simple 3-manifold, and let $Y$ and 
$Y'$ be distinct boundary components of
$M$. Suppose that there is an orientation-preserving
homeomorphism $h \colon Y \rightarrow S$ and an
orientation-reversing homeomorphism $h' \colon
S \rightarrow Y'$, where $S$ is some surface of
genus at least two. Let $\psi \colon S \rightarrow S$
be a pseudo-Anosov homeomorphism, and let
$M/\!\sim$ be the manifold obtained by
gluing $Y$ and $Y'$
via $h' \psi^n h$. Then,
provided $|n|$ is sufficiently large,
$$g(M) - g(S) + 1 \leq g(M/\!\sim) \leq g(M) + g(S) + 1.$$

}

The proof is very similar, but not identical, to
that of the main theorem. To achieve the
upper bound on $g(M/\!\sim)$, one starts with
a minimal genus splitting for $M$, and uses it
to construct a splitting for $M/\!\sim$.
One might have to modify the surface in $M$
to ensure that it does not separate $Y$ from $Y'$.
This may increase its genus by $g(S)$. Then,
to construct a Heegaard surface for $M/\!\sim$,
one attaches a tube that runs through $S$.
This increases the genus of the surface by one. Hence, we obtain
the upper bound. An instructive example
is where $M$ is the product $S \times I$ of a closed orientable
surface and an interval, and where $M/\!\sim$ 
fibres over the circle. (Of course,
though, $M$ is not simple in this case.)
Then, $g(M) = g(S)$, but in general,
$g(M/\!\sim)$ may be as much as $2g(S) + 1$. (See [6] for example).

To achieve the lower bound on $g(M/\!\sim)$, one starts with
a minimal genus Heegaard surface $F$ for $M/\!\sim$.
One untelescopes it to a generalised Heegaard
splitting satisfying the five conditions given earlier.
Using the geometry of $M/\!\sim$, one
can show that this is disjoint from a copy of
$S$ in $M/\!\sim$, provided $|n|$ is sufficiently large.
Thus, it gives a generalised
Heegaard splitting for $M$, which can be amalgamated
to form a Heegaard surface. One calculates its genus
to be $g(F) + g(S) - 1$.

The same issues arise when gluing simple manifolds $M$
and $M'$ but when $\partial M$ and $\partial M'$
are disconnected. Again, one does not obtain
an exact equality.

It should be possible to generalise the main theorem
even further. One can consider the manifold
$M \cup_{h' \psi h} M'$, where $\psi \colon S \rightarrow
S$ is some homeomorphism. It should be true that,
under the hypotheses of the main theorem,
and provided the distance of $\psi$ is sufficiently
large, then the conclusion of the main theorem holds.
Here, distance is as measured by the action of
$\psi$ on the curve complex of $S$. This would
indeed represent a generalisation, since the
distance of $\psi^n$, for a given pseudo-Anosov
$\psi$, is arbitrarily large, provided $|n|$
is sufficiently large [2]. One might also try
to drop the assumptions that $M$ and $M'$
are acylindrical, or even that they have
incompressible boundary. But one would then
need to make further hypotheses on $\psi$.
To prove these more general results, one would need 
to establish geometric control on $M \cup M'$, using 
the theory of Kleinian groups.

\vskip 18pt
\centerline {\caps References}
\vskip 6pt

\item{1.} {\caps M. Freedman, J. Hass and P. Scott,}
{\sl Least area incompressible surfaces in $3$-manifolds},
Invent. Math. {\bf 71} (1983) 609--642.

\item{2.} {\caps J. Hempel}, {\sl 3-manifolds as viewed
from the curve complex}, Topology {\bf 40} (2001)
631--657.

\item{3.} {\caps K. Johannson}, {\sl Topology and combinatorics of
3-manifolds}, Springer-Verlag (1995).

\item{4.} {\caps M. Lackenby}, {\sl Heegaard splittings,
the virtually Haken conjecture and Property $(\tau)$},
Preprint.

\item{5.} {\caps J. Pitts and J. H. Rubinstein,} {\sl Existence 
of minimal surfaces of bounded topological type in three-manifolds.}
Miniconference on geometry and partial differential equations 
(Canberra, 1985), 163--176.

\item{6.} {\caps J. H. Rubinstein}, {\sl Minimal surfaces
in geometric 3-manifolds}, Preprint.

\item{7.} {\caps M. Scharlemann,} {\sl Heegaard splittings
of compact 3-manifolds}, Handbook of Geometric Topology,
Elsevier (2002), 921--953.

\item{8.} {\caps M. Scharlemann and A. Thompson}, {\sl
Thin position for $3$-manifolds,} Geometric topology (Haifa,
1992), 231--238, Contemp. Math., 164.

\item{9.} {\caps R. Schoen and S. Yau}, {\sl Existence of
incompressible surfaces and the topology of 3-manifolds
with non-negative scalar curvature}, Ann. Math. {\bf 110}
(1979) 127--142.

\item{10.} {\caps J. Schultens}, {\sl The classification of
Heegaard splittings for (compact orientable surface)$\times S^1$},
Proc. London Math. Soc. {\bf 67} (1993), 425--448.

\item{11.} {\caps J. Schultens}, {\sl Heegaard genus formula for
Haken manifolds}, Preprint.

\item{12.} {\caps T. Soma}, {\sl Volume of hyperbolic 3-manifolds
with iterated pseudo-Anosov amalgamations}, Geom. Ded. {\bf 90}
(2002), 183--200.

\vskip 12pt
\+ Mathematical Institute, Oxford University, \cr
\+ 24-29 St Giles', Oxford OX1 3LB, United Kingdom. \cr

\end